\documentclass{svjour3}     
\usepackage{placeins}
\usepackage{graphicx}
\usepackage{amsmath}
\usepackage{epstopdf} 
\usepackage[ansinew]{inputenc}
\usepackage{amssymb}
\usepackage{makeidx}
\usepackage{float}
\restylefloat{table}

\usepackage{color}
\usepackage{listings}
\definecolor{juliaBackground}{RGB}{255,250,255}
\definecolor{juliaComment}{RGB}{155,142,226}
\definecolor{juliaString}{RGB}{215,0,122}
\definecolor{juliakeyword}{RGB}{158,43,56}

\lstdefinelanguage{Julia}
  {morekeywords={abstract,break,case,catch,const,continue,do,else,elseif,
      end,export,false,for,function,immutable,import,importall,if,in,
      macro,module,otherwise,quote,return,switch,true,try,type,typealias,
      using,while,range,plot,display,fit,exp},
   sensitive=true,
   morecomment=[l]\#,
   morecomment=[n]{\#=}{=\#},
   morestring=[s]{"}{"},
   morestring=[s]{'}{'},
}[keywords,comments,strings]

\lstdefinestyle{Julia}{
    language=Julia,
    backgroundcolor=\color{juliaBackground},
    commentstyle=\color{juliaComment},
    stringstyle =\color{juliaString},
    keywordstyle = \color{juliakeyword},
}
\usepackage[autoload=false]{jlcode}

\title{Construction of DAE solver Rodas4PW and numerical test examples with Rosenbrock and W-methods in OrdinaryDiffEqRosenbrock.jl}
\titlerunning{Construction of DAE solver Rodas4PW} 

\author{Gerd Steinebach}
\date{\today}

\begin{document}
\maketitle

\begin{abstract}
This paper compares various Rosenbrock and Rosenbrock W-methods. To achieve this, some methods implemented in the Julia package OrdinaryDiffEqRosenbrock.jl 
are selected and benchmarked using test examples. The test examples include stiff ODE problems as well as differential-algebraic equations.
In addition to the already known methods, the new W-method {\tt Rodas4PW} is introduced and included in the tests. 
\end{abstract}

\section{Rosenbrock W-methods}
Rosenbrock, also known as Rosenbrock-Wanner (ROW) methods for solving stiff ordinary differential equations (ODEs) and differential-algebraic equations (DAEs) 
have been well researched, and a large number of efficient methods is available.
The theory behind these methods is described in great detail in Hairer \& Wanner \cite{hairer}, and a good overview of the development of 
the Rosenbrock methods can be found in Lang \cite{lang}.

Whilst the early methods, such as GRK4T \cite{kaps}, were only suitable for ODEs, 
since Roche \cite{roche} most Rosenbrock methods have 
also been designed for DAEs of the following mass matrix form:
\begin{equation}
M y' = f(t,y) \; , \quad y(t_0)=y_0 \; . \label{eq:dae} 
\end{equation}  
In order to solve DAEs, the methods must satisfy additional order conditions.
Usually, the DAE problems must be of index 1, see \cite{hairer}.
A Rosenbrock method for solving (\ref{eq:dae}) is given by:
\begin{align}
(M - h \, \gamma \, f_y) k_i &= h \, f(t_0+ \alpha_i h, y_0 + \sum_{j=1}^{i-1} \alpha_{ij} k_j) 
+ h \, f_y \sum_{j=1}^{i-1} \gamma_{ij}  k_j + h^2 \gamma_i f_t , \label{eq_row1} 
\end{align}
\begin{align}
y_1 &= y_0 + \sum_{i=1}^{s} b_i k_i,  \label{eq_row2} \quad 
\mbox{with} \quad f_y = \frac{\partial f}{\partial y}(t_0,y_0) \; , \quad f_t = \frac{\partial f}{\partial t}(t_0,y_0) \, . \nonumber
\end{align}
$y_1$ defines the approximation of the solution $y(t_0+h)$ with step size $h$.
The coefficients and weights of the method are $\gamma$, $\alpha_{ij}$, $\gamma_{ij}$, and $b_i$. 
Moreover, $\alpha_i = \sum_{j=1}^{i-1} \alpha_{ij}$ and
$\gamma_i = \gamma + \sum_{j=1}^{i-1} \gamma_{ij}$ holds.\\
W-methods \cite{steihaug} overcome the drawback of Rosenbrock methods, which require the Jacobian matrix to be recalculated at every time step.
They can use an approximation $J \approx f_y$ of the Jacobian.
In the simplest case, this is achieved by keeping the matrix constant over several time steps. For stability reasons, however, this number should not be set too high.
In particular, for DAE problems, $J$ cannot be chosen arbitrarily. Here, at least the derivatives of the algebraic 
equations with respect to the algebraic state variables must be exact.
To achieve this, the system (\ref{eq:dae}) must be split into differential and algebraic components:
\begin{eqnarray}
y_1' &=& f_1(t,y_1,y_2) \, , \label{eq:dae1} \\
0 &=& f_2(t,y_1,y_2) \, . \label{eq:dae2} 
\end{eqnarray}
When $M$ is not a diagonal matrix this can be done using a singular value decomposition of $M$. 

The disadvantage of W-methods is that the coefficients of the methods must satisfy significantly more order conditions.
In his thesis, Jax \cite{jax2} derives the necessary conditions. For example, for an order-three method, a total 
of 26 order conditions must be satisfied. By contrast, a Rosenbrock method with an exact Jacobian matrix requires 
only 6 order conditions.

The following Table \ref{tab:overview} is updated from \cite{steijulia} and lists all Rosenbrock-type methods currently 
implemented in the package {\tt OrdinaryDiffEqRosenbrock.jl}. 
In addition to the number of stages $s$ for the methods, their order of convergence for ODEs and DAEs is given. 
When the methods are applied to semi-discretised parabolic partial differential equations and to the Prothero–Robinson model, 
an order reduction usually occurs \cite{oro,rang2,scholz}. The order achieved for a parabolic model problem is therefore also listed. 
Furthermore, the W-property is specified, i.e. the order achieved for ODEs or DAEs when the Jacobian matrix is not exact. 
The values for DAEs given here assume that only the derivatives  
$\frac{\partial f_2}{\partial y_2}$ of the algebraic equations (\ref{eq:dae2}) are exact.
It is also specified whether the methods are A-stable and L-stable, i.e. $R(\infty)=0$, where $R$ is the stability function.
All values given in brackets refer to the properties of the embedded method used for step size control. 
All methods belonging to the {\tt Rodas} family, as well as {\tt Rosenbrock23}, {\tt Rosenbrock32}, are equipped 
with their own interpolation, which allows for a continuous, dense output of the solution.

\begin{table*}[t] 
\caption{Rosenbrock-type methods within \texttt{OrdinaryDiffEqRosenbrock.jl} and their properties.}
\begin{tabular}{r|ccccccccc}
Method & Ref. & $s$ & ODE & W & DAE & Para & A& $R(\infty)$ & Cont \\
\hline
GRK4T &\cite{kaps} 1979&4&4(3)&-&2(-)&2(-)&-&0.45(2.6)&- \\
GRK4A &\cite{kaps} 1979&4&4(3)&-&2(2)&2(2)&x&0.99(0.31)&- \\
RosShamp4 &\cite{shamp} 1982&4&4(3)&-&2(2)&2(2)&x&0.33(1.0)&- \\
Veldd4 &\cite{veld} 1984&4&4(3)&-&2(-)&2(-)&-&0.24(2.6)& -\\
Velds4 &\cite{veld} 1984&4&4(3)&2/-&2(2)&2(2)&x&0.33(0.33)& -\\
Scholz4\_7 &\cite{scholz} 1989&3&3(1)&-&3(2)&3(2)&x&0.73(0.52)&- \\
Ros4LStab &\cite{hairer} 1991&4&4(3)&-&2(2)&2(2)&x&0.0(0.55)& -\\
Rodas4 &\cite{hairer} 1991&6&4(3)&-&4(3)&2(2)&x&0.0(0.0)&x \\
Rodas42 &\cite{hairer} 1991&6&4(3)&-&4(3)&2(2)&x&0.0(0.0)&x \\
Rodas5 &\cite{rodas5} 1993&8&5(4)&-&5(4)&2(2)& x&0.0(0.0)&x \\
Rodas4P &\cite{rodas4p} 1995&6&4(3)&-&4(3)&4(3)&x&0.0(0.0)&x \\
ROS3 &\cite{rodas3} 1997&3&3(2)&-&2(2)&2(2)&x&0.0(0.5)& -\\
Rodas3 &\cite{rodas3} 1997&4&3(2)&-&3(2)&2(2)&x&0.0(0.0)&- \\
Rosenbrock23 &\cite{shampine} 1997&3&2(3)&2/-&2(-)&2(-)&x&0.0(1.6)& x\\
Rosenbrock32 &\cite{shampine} 1997&3&3(2)&2/-&-(2)&-(2)&-&1.6(0.0)& x\\
ROS3P &\cite{lang2} 2001&3&3(2)&-&3(2)&3(2)&x&0.73(0.73)&- \\
ROS34PW1a &\cite{ranga} 2005&4&3(2)&3/-&3(2)&3(2)&x&0.0(0.0)&- \\
ROS34PW1b &\cite{ranga} 2005&4&3(2)&3/-&3(2)&3(2)&x&0.0(0.0)&- \\
ROS34PW2 &\cite{ranga} 2005&4&3(2)&3/2&3(2)&3(2)&x&0.0(0.48)&- \\ 
ROS34PW3 &\cite{ranga} 2005&4&4(2)&3/-&3(2)&3(2)&x&0.63(0.43)&- \\
ROS2S &\cite{rang12} 2012&3&2(1)&2/-&2(2)&2(1)&x&0.0(0.33)&- \\
ROS2PR &\cite{rang12a} 2012&3&2(1)&-&2(2)&2(1)&x&0.0(0.0)&- \\
ROS34PRw &\cite{rang12a} 2012&4&3(2)&3/2&3(2)&3(2)&x&0.0(0.25)&- \\ 
ROS3PR &\cite{rang1} 2014 &3& 3(2)&-&3(2)&3(2)&x&0.73(0.73)& -\\
ROS3PRL &\cite{rang1} 2014&4&3(2)&-&3(2)&3(2)&x&0.0(0.25)&- \\
ROS3PRL2 &\cite{rang15} 2015&4&3(2)&-&3(2)&3(2)&x&0.0(0.25)&- \\
Rodas4P2 &\cite{rodas4p2} 2020&6&4(3)&2/2&4(3)&4(3)&x&0.0(0.0)&x \\ 
Rodas5P &\cite{rodas5p} 2023&8&5(4)&2/2&5(4)&4(4)& x&0.0(0.0)&x \\ 
Rodas23W &\cite{steijulia} 2024&5&2(3)&2/2&2(3)&2(3)&x&0.0(0.0)& x\\
Rodas3P  &\cite{steijulia} 2024&5&3(2)& - &3(2)&3(2)&x&0.0(0.0)& x\\
Rodas5Pe &\cite{steijulia} 2024&8&5(4)&2/2&5(3)&4(4)& x&0.0(0.46)&x \\ 
Rodas5Pr &\cite{steijulia} 2024&8&5(4)&2/2&5(4)&4(4)& x&0.0(0.0)&x \\ 
Rodas3d &\cite{rodas3d} 2024& 4&3(2)&-&3(2)&2(2) &x& 0.0(0.0) & - \\ 
Rodas6P &\cite{rodas6p} 2025&19&6(5)&2/2&6(5)&4(3)&x&0.0(0.0)&x \\ 
Tsit5DA &\cite{rodas6p} 2025&12&5(4)&-&5(4)&-&-&-&x \\ 
\hline 
Rodas4PW &2026&9&4(3)&4/2&4(3)&4(3)&x&0.0(0.0)&x \\ 
\hline 
\end{tabular}
\label{tab:overview}

s = number of stages, order of convergence for ODEs, DAEs and parabolic problem (embedded
method in brackets), W = order of W-method for ODEs / DAEs, A = A stable, $R(\infty)$ = limit of stability function, Cont = continuous output for DAEs.
\end{table*}

The outline of this paper is as follows: The next section introduces the new method {\tt Rodas4PW} and discusses 
the implementation of W-methods. In Section 3, a number of benchmarks are presented and selected methods are tested 
against them. In particular, the aim is to identify when the use of W-methods pays off compared to ROW methods. 
Complete Julia programmes are provided for all benchmarks, so that the reader can carry out own tests using 
other methods.

\section{Order conditions, construction and implementation of Rodas4PW}

The aim is to construct a fourth-order W-method for DAE problems. 
The corresponding order conditions are listed in Table \ref{tab:rodas4pw} together with their Butcher trees. 
These were derived in \cite{jaxstei} and guarantee order $p=4$ for problem 
(\ref{eq:dae1}),(\ref{eq:dae2}), when the components $\frac{\partial f_2}{\partial y_1}$, $\frac{\partial f_2}{\partial y_2}$ of the Jacobian
are assumed to be exact, while $\frac{\partial f_1}{\partial y_1}$, $\frac{\partial f_1}{\partial y_2}$ may be arbitrary. 
If, additionally $\frac{\partial f_2}{\partial y_1}$ is not exact, the order is reduced to $p=2$.
According to Jax \cite{jax2}, W-methods that allow for approximations of  $\frac{\partial f_2}{\partial y_1}$ are also possible. 
However, the number of order conditions to be satisfied then increases very sharply, and no methods of order higher 
than $p=3$ are currently known.

{\tt Rodas4PW}, when used as a ROW method with an exact Jacobian matrix at each time step, should furthermore preserve the properties of {\tt Rodas4P2}
and {\tt Rodas5P}, respectively. These are stiffly accurate \cite{hairer}
and the order reduction for parabolic problems is avoided. 
Therefore, the additional order conditions of Scholz \cite{scholz} are satisfied, see \cite{rodas4p2}. 
In total, 66 non-linear equations had to be satisfied for the method, including the embedded scheme of order $\hat p = 3$. 
Chosing a stage number of $s=9$ a solution to the system of all order conditions was found using the methods of the {\tt NonlinearSolve.jl} package.
In addition, {\tt Rodas4PW}, like all the other methods in the Rodas family, allows
for dense output of the solution via its own interpolation:
\begin{eqnarray}
 y(t_0 + \tau \, h) &=& y_0 + \sum_{i=1}^s b_i(\tau) k_i \; , \quad \tau \in [0,1] \, ,\label{eq:interp1} \\
 b_i(\tau) &=& \tau (b_i-c_i) + \tau^2 (c_i-d_i) + \tau^3 d_i \, .\label{eq:interp2}
\end{eqnarray}
All coefficients of {\tt Rodas4PW} are given in Listing \ref{lst:coeffs} in the Appendix.

To make use of the W property of the method, a wide range of implementations are possible. 
The specific implementation depends heavily on the problem to be solved. Only a small selection of the possibilities 
is discussed here. There is no universal answer to the question of how W-methods can best achieve efficiency gains over ROW methods. 

\begin{enumerate}
\item Freezing the Jacobian matrix over several time steps. This is also the implementation chosen in 
{\tt OrdinaryDiffEqRosenbrock.jl} for W-methods when applied to ODEs. In this approach, the Jacobian matrix is kept constant 
for as long as the current timestep $h$ satisfies $|\frac{h}{h_{last}} - 1| > \mbox{\tt jac\_reuse\_gamma\_tol}$. 
In addition, a maximum number {\tt max\_jac\_age} of time steps for 
the freezing is defined and, whenever a timestep is rejected, a new Jacobian is computed.
The current default settings are {\tt jac\_reuse\_gamma\_tol} = 0.03, {\tt max\_jac\_age} = 20.
\item
Freezing the W-matrix.
In ODE problems, the matrix $W=(I - h \, \gamma \, f_y)$ is kept constant over multiple time steps. This eliminates the need to compute a new LU decomposition 
of W for the systems of linear equations to be solved in (\ref{eq_row1}) at every time step.
\item Simplification of the structure of the Jacobian matrix.
By omitting non-zero entries, the matrix - and thus also $W$ can be reduced, for example, to a banded form. 
For DAEs of the form (\ref{eq:dae1},\ref{eq:dae2}), the W matrix is: 
\[ W = \left( \begin{array}{cc} I - h \gamma \frac{\partial f_1}{\partial y_1} & -h \gamma \frac{\partial f_1}{\partial y_2} \\
-h \gamma \frac{\partial f_2}{\partial y_1} & -h \gamma \frac{\partial f_2}{\partial y_2} \end{array} \right) \]
If $\frac{\partial f_1}{\partial y_2}$ is set to zero, the block-linear system of equations with $W$ can be easily reduced to two linear systems of smaller dimension.
{\tt Tsit5DA} is a method in which $\frac{\partial f_1}{\partial y_1}$ and $\frac{\partial f_1}{\partial y_2}$ are set both to zero, \cite{rodas6p}.
\item Operator splitting. 
For ODE problems of the form $y' = f(y) + g(y)$, it may be useful to treat the $f$ and $g$ terms separately. For example, in a convection-reaction problem, 
one would compute the convection term $f$ explicitly and the reaction term $g$ implicitly. Such an IMEX (implicit-explicit) approach can be implemented as the 
W-method based on the Jacobian matrix by setting $J = \frac{\partial g}{\partial y}$.
\item Approximate matrix factorization (AMF).
For example, in a diffusion-reaction problem, the Jacobian matrix $J = J_f + J_g$ can be decomposed into a constant linear operator 
$J_f=\frac{\partial f}{\partial y}$ and the stiff reaction term $J_g=\frac{\partial g}{\partial y}$.
The matrix $W=(I - h \gamma J)$ is then approximated by $W = (I - h \gamma J_f)(I-h \gamma J_g)$ and the two resulting systems of linear equations can be solved separately. 
Another option would be, for example, to decompose $J = J_L + J_U$ additively into a lower and an upper triangular matrix and set 
$W = (I - h \gamma J_L)(I-h \gamma J_R)$.\\
Or, for example by using a singular value decomposition, one can obtain a low-rank approximation of $J$ given by $J \approx A\, B \,C$, 
where $A$, $B$, $C$ are matrices of dimension $n \times s$, $s \times s$, $s \times n$, respectively. 
Since
\[ (I - h \gamma A \, B \, C)^{-1} = I + h \gamma A (I - h \gamma B \, C \, A)^{-1} B \, C \]
only linear systems of dimension $s<n$ have to be solved.
\end{enumerate}
This overview of implementation options is certainly not exhaustive; see also \cite{jax2}. 
No systematic investigations are conducted here, instead, only some numerical observations are discussed in the next section.

\section{Benchmarks} 

First, we compute simple order tests taken from \cite{rodas5p,rodas6p}. We solve the DAE
\begin{equation}
y_1' = \frac{y_2}{y_1} \; ,\quad
0 = \frac{y_1}{y_2} - t \, , \label{eq:prob1}
\end{equation}
with exact solution $y_1(t)=\ln(t)$, $y_2(t)=\frac{\ln(t)}{t}$ in the time interval $t \in [2,4]$.
Methods {\tt Rodas4P2} \cite{rodas4p2}, {\tt Rodas5P} \cite{rodas5p}, {\tt Rodas6P} \cite{rodas6p}, 
{\tt ROS34PW2} \cite{ranga}, {\tt ROS34PRw} \cite{rang12a} and the new method {\tt Rodas4PW} are applied to this problem with different constant step sizes. 
Listing \ref{lst:order_test} shows the implementation and part a) of Table \ref{tab:prob1} the obtained absolute errors and the numerically orders at final time $t=4$. 

\FloatBarrier
\begin{lstlisting}[caption=Order tests., label=lst:order_test, style=Julia]
using OrdinaryDiffEqRosenbrock, PrettyTables

function fcn!(dy, y, p, t)
  dy[1] = y[2]/y[1]
  dy[2] = y[1]/y[2] - t
end

function anasol(t)
  [log(t), log(t)/t]
end

tspan = (2.0, 4.0); u0 = anasol(tspan[1]); p = nothing
M = zeros(2,2); M[1,1] = 1;
prob = ODEProblem(ODEFunction(fcn!, mass_matrix=M), u0, tspan, p) 

N = [8,16,32,64,128,256]
Algs = [Rodas4P2(), Rodas5P(), Rodas6P(), ROS34PW2(), ROS34PRw(), Rodas4PW()]  
Algs_string = [nameof(typeof(a)) for a in Algs]

N_N = length(N); N_Algs = length(Algs)
errors = zeros(N_N, N_Algs); orders = zeros(N_N - 1, N_Algs)

for ialg = 1 : N_Algs
  for n = 1 : N_N
    dt = (tspan[2] - tspan[1])/N[n]
    sol = solve(prob, Algs[ialg], adaptive = false, dt = dt) 
    errors[n, ialg] = maximum(abs.(anasol(tspan[2]) - sol(tspan[2])))
    if n > 1  
        orders[n-1, ialg] = log(errors[n-1, ialg] / errors[n, ialg]) / log(2) 
    end
  end
end
pretty_table(errors, column_labels = Algs_string, formatters = [fmt__printf("%.2e")], display_size= (40,188))
pretty_table(orders, column_labels = Algs_string, formatters = [fmt__printf("%.2f")], display_size= (40,188))
\end{lstlisting}


\begin{table}
\caption{Numerical errors and orders of the various methods for problem (\ref{eq:prob1}).}
\label{tab:prob1}
\begin{scriptsize}

a) Updated Jacobian at each time step

\begin{center}
\begin{tabular}{|r|r|r|r|r|r|}
  \hline
  \textbf{Rodas4P2} & \textbf{Rodas5P} & \textbf{Rodas6P} & \textbf{ROS34PW2} & \textbf{ROS34PRw} & \textbf{Rodas4PW} \\
  \hline
  4.16e-06 & 1.07e-06 & 3.80e-08 & 3.66e-05 & 4.75e-05 & 4.07e-06 \\
  2.20e-07 & 2.93e-08 & 5.03e-10 & 3.79e-06 & 7.81e-06 & 2.09e-07 \\
  1.29e-08 & 8.56e-10 & 7.25e-12 & 5.21e-07 & 1.09e-06 & 1.54e-08 \\
  7.81e-10 & 2.59e-11 & 1.08e-13 & 6.75e-08 & 1.44e-07 & 1.05e-09 \\
  4.82e-11 & 8.03e-13 & 2.22e-15 & 8.57e-09 & 1.84e-08 & 6.88e-11 \\
  2.99e-12 & 2.69e-14 & 2.44e-15 & 1.08e-09 & 2.33e-09 & 4.60e-12 \\
  \hline
\end{tabular}
\begin{tabular}{|r|r|r|r|r|r|}
  \hline
  \textbf{Rodas4P2} & \textbf{Rodas5P} & \textbf{Rodas6P} & \textbf{ROS34PW2} & \textbf{ROS34PRw} & \textbf{Rodas4PW} \\
  \hline
  4.24 & 5.19 & 6.24 & 3.27 & 2.60 & 4.29 \\
  4.10 & 5.10 & 6.11 & 2.86 & 2.84 & 3.76 \\
  4.04 & 5.05 & 6.07 & 2.95 & 2.93 & 3.87 \\
  4.02 & 5.01 & 5.60 & 2.98 & 2.96 & 3.94 \\
  4.01 & 4.90 & -0.14 & 2.99 & 2.98 & 3.90 \\
  \hline
\end{tabular}
\end{center}

b) Differential parts $\frac{\partial f_1}{\partial y_1}$, $\frac{\partial f_1}{\partial y_2}$ are constant

\begin{center}
\begin{tabular}{|r|r|r|r|r|r|}
  \hline
  \textbf{Rodas4P2} & \textbf{Rodas5P} & \textbf{Rodas6P} & \textbf{ROS34PW2} & \textbf{ROS34PRw} & \textbf{Rodas4PW} \\
  \hline
  3.47e-05 & 5.17e-05 & 1.80e-05 & 1.67e-04 & 7.82e-05 & 1.66e-05 \\
  2.47e-05 & 2.02e-05 & 5.89e-06 & 1.73e-05 & 8.38e-06 & 9.67e-07 \\
  8.07e-06 & 5.97e-06 & 1.66e-06 & 1.95e-06 & 9.94e-07 & 5.79e-08 \\
  2.25e-06 & 1.61e-06 & 4.40e-07 & 2.31e-07 & 1.22e-07 & 3.54e-09 \\
  5.91e-07 & 4.16e-07 & 1.13e-07 & 2.81e-08 & 1.51e-08 & 2.19e-10 \\
  1.51e-07 & 1.06e-07 & 2.87e-08 & 3.46e-09 & 1.88e-09 & 1.35e-11 \\
  \hline
\end{tabular}
\begin{tabular}{|r|r|r|r|r|r|}
  \hline
  \textbf{Rodas4P2} & \textbf{Rodas5P} & \textbf{Rodas6P} & \textbf{ROS34PW2} & \textbf{ROS34PRw} & \textbf{Rodas4PW} \\
  \hline
  0.49 & 1.35 & 1.61 & 3.27 & 3.22 & 4.10 \\
  1.61 & 1.76 & 1.82 & 3.15 & 3.08 & 4.06 \\
  1.84 & 1.89 & 1.92 & 3.08 & 3.03 & 4.03 \\
  1.93 & 1.95 & 1.96 & 3.04 & 3.01 & 4.02 \\
  1.97 & 1.98 & 1.98 & 3.02 & 3.01 & 4.01 \\
  \hline
\end{tabular}
\end{center}

c) Differential parts $\frac{\partial f_1}{\partial y_1}$, $\frac{\partial f_1}{\partial y_2}$ and algebraic part $\frac{\partial f_2}{\partial y_1}$ are constant

\begin{center}
\begin{tabular}{|r|r|r|r|r|r|}
  \hline
  \textbf{Rodas4P2} & \textbf{Rodas5P} & \textbf{Rodas6P} & \textbf{ROS34PW2} & \textbf{ROS34PRw} & \textbf{Rodas4PW} \\
  \hline
  7.63e-05 & 9.49e-06 & 1.73e-05 & 1.15e-04 & 6.97e-05 & 2.83e-05 \\
  3.29e-05 & 9.50e-06 & 5.19e-06 & 3.19e-05 & 1.51e-05 & 8.17e-06 \\
  9.86e-06 & 3.27e-06 & 1.42e-06 & 1.02e-05 & 4.28e-06 & 2.13e-06 \\
  2.67e-06 & 9.28e-07 & 3.72e-07 & 2.76e-06 & 1.11e-06 & 5.41e-07 \\
  6.91e-07 & 2.46e-07 & 9.51e-08 & 7.13e-07 & 2.82e-07 & 1.36e-07 \\
  1.76e-07 & 6.32e-08 & 2.41e-08 & 1.81e-07 & 7.10e-08 & 3.42e-08 \\
  \hline
\end{tabular}
\begin{tabular}{|r|r|r|r|r|r|}
  \hline
  \textbf{Rodas4P2} & \textbf{Rodas5P} & \textbf{Rodas6P} & \textbf{ROS34PW2} & \textbf{ROS34PRw} & \textbf{Rodas4PW} \\
  \hline
  1.21 & -0.00 & 1.74 & 1.85 & 2.20 & 1.79 \\
  1.74 & 1.54 & 1.87 & 1.65 & 1.82 & 1.94 \\
  1.89 & 1.82 & 1.93 & 1.88 & 1.94 & 1.98 \\
  1.95 & 1.92 & 1.97 & 1.95 & 1.98 & 1.99 \\
  1.97 & 1.96 & 1.98 & 1.98 & 1.99 & 2.00 \\
  \hline
\end{tabular}
\end{center}
\end{scriptsize}
\end{table}

\FloatBarrier
\begin{lstlisting}[caption=Order tests with constant differential components of the Jacobian matrix., label=lst:order_test_w, style=Julia]
function jac!(J, y, p, t)
    y0 = anasol(2.0) 
    J[1, 1] = - y0[2] / y0[1]^2 
    J[1, 2] = 1 / y0[1]
    J[2, 1] = 1 / y[2]
    J[2, 2] = -y[1] / y[2]^2
end

prob = ODEProblem(ODEFunction(fcn!, mass_matrix=M, jac=jac!), u0, tspan, p)
\end{lstlisting}
Next, we freeze the differential parts of the Jacobian. The modified implementation is given in Listing \ref{lst:order_test_w} and the results in part b) of 
Table \ref{tab:prob1}. Finally, additionally the algebraic part $\frac{\partial f_2}{\partial y_1}$ of the Jacobian is held constant.
We can observe in Table \ref{tab:prob1} that all theoretical orders are being achieved.

The next example is taken from the documentation of {\tt OrdinaryDiffEq.jl}.
The two dimensional diffusion-reaction equation reads
 \begin{equation}
\frac{\partial u}{\partial t} = D \left( \frac{\partial^2 u}{\partial x^2} + \frac{\partial^2 u}{\partial y^2}\right) + u^2(1-u) + e^x\, , \label{eq:prob2}
\end{equation}
with initial values $u(x,y,0) = 16 \, x \, y \,(1-x)(1-y)$ and vanishing boundary conditions on the domain $(x,y) \in [0,1]^2$.
We choose $D=0.1$ and the time interval $t \in [0,1]$.

\begin{figure}[t]
 \centering
 \includegraphics[width=0.49\textwidth]{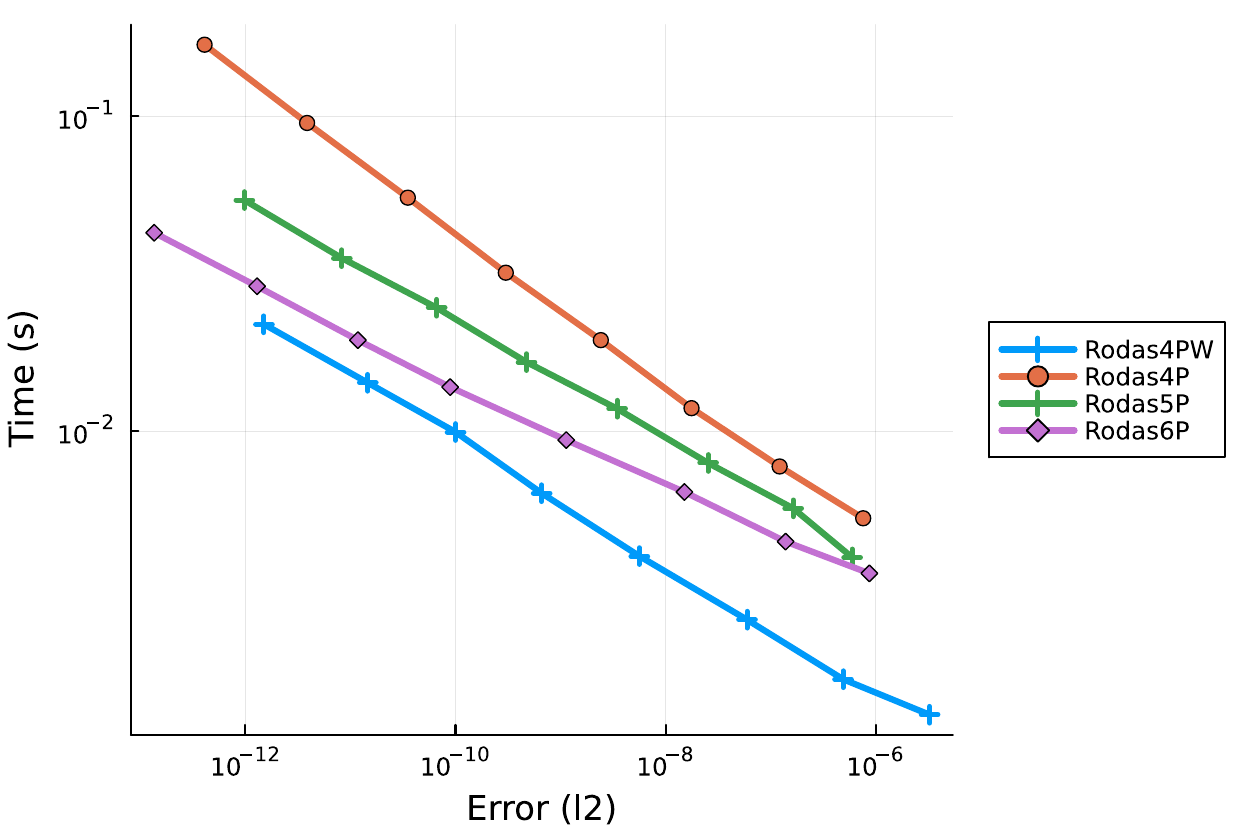} \hfill
 \includegraphics[width=0.49\textwidth]{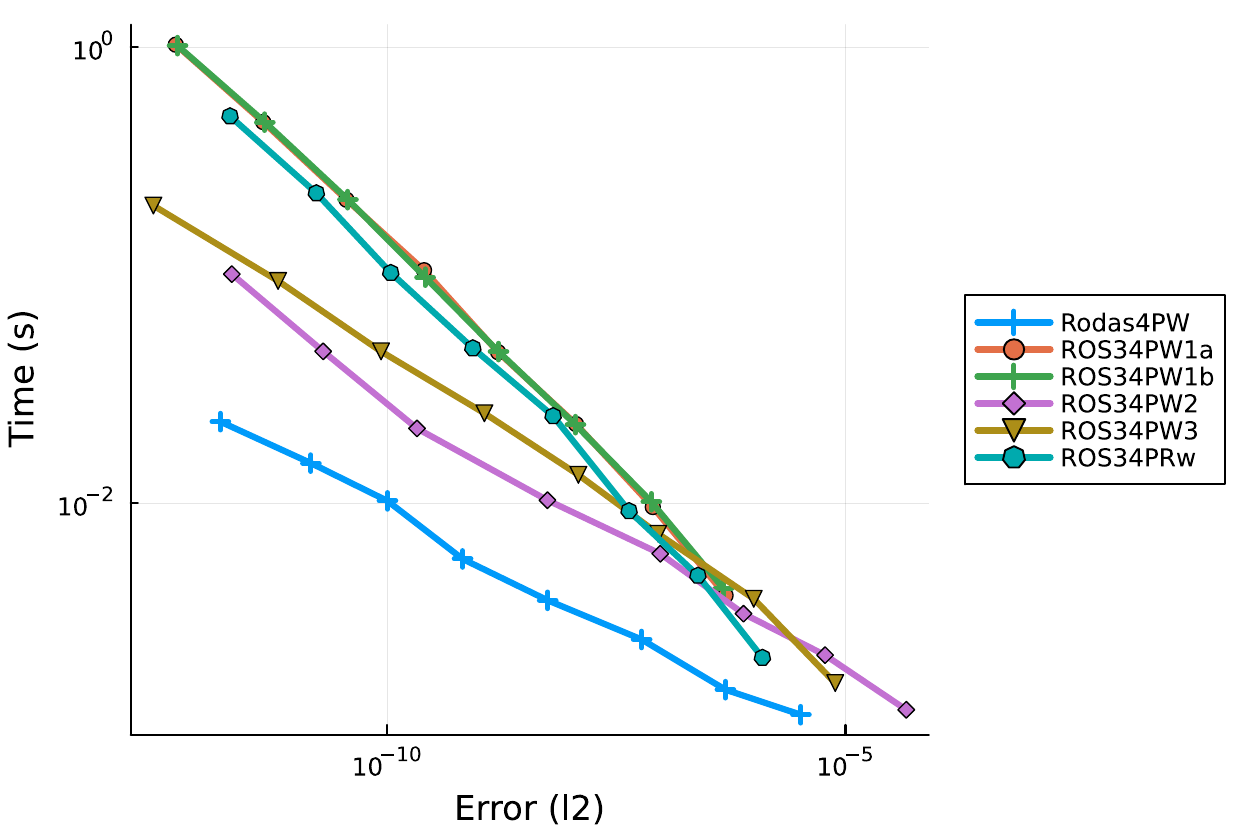}
 \caption{Results of two dimensional diffusion-reaction benchmark with inefficient implementation.}\label{fig:bench_1}
\end{figure}

\begin{figure}[t]
 \centering
 \includegraphics[width=0.65\textwidth]{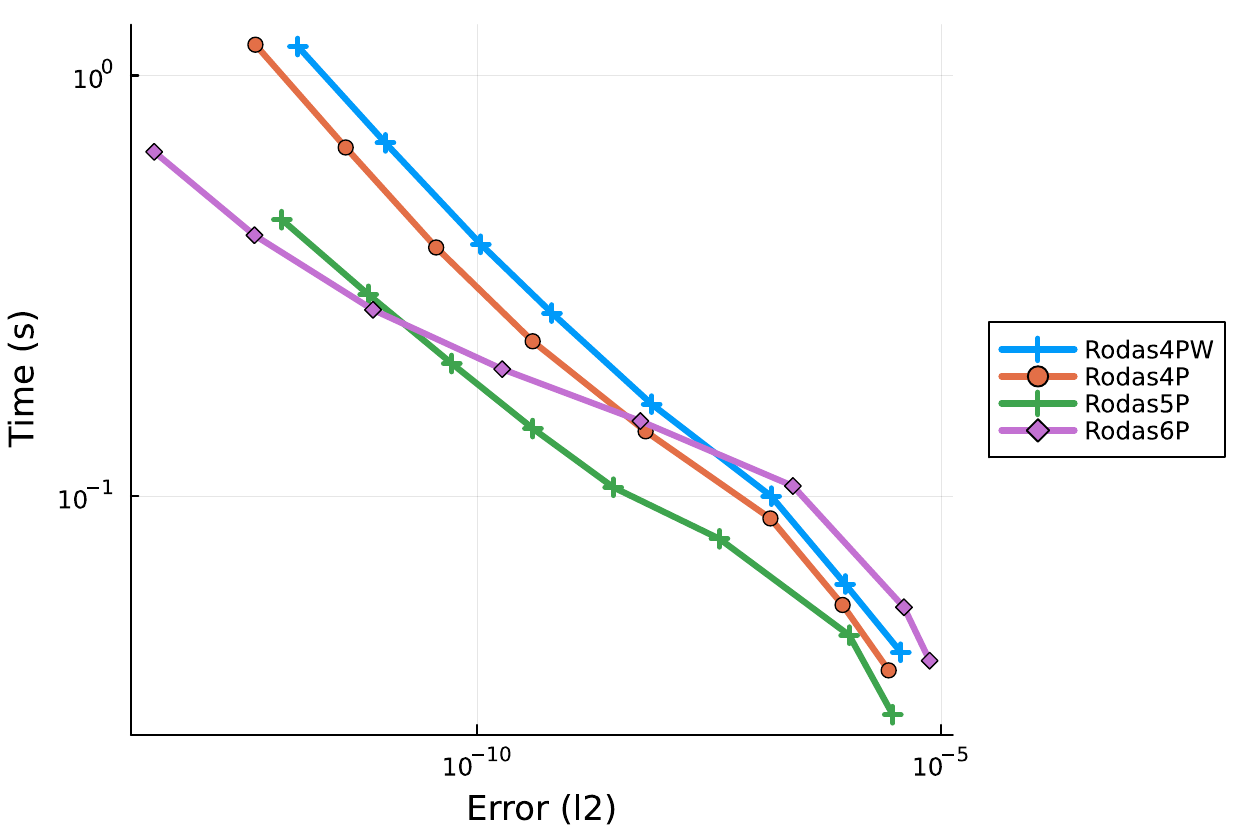} 
 \caption{Results of two dimensional diffusion-reaction benchmark with sparse matrix operations.}\label{fig:bench_2}
\end{figure}

We run two work-precision benchmarks: one using a naive, inefficient implementation, and one using an implementation that takes advantage of the sparse Jacobian matrix.
In both cases, the space derivatives are approximated using central differences. The first implementation, using a 10x10 space grid, is 
shown in Listing \ref{lst:bench_1}, and the results are displayed in Figure \ref{fig:bench_1}.
The second implementation uses a 60x60 grid and defines the diffusion matrix as sparse. The structure of the matrix is communicated to the 
integrator as the pattern of the Jacobian, and the KLU decomposition -designed for sparse matrices- is used as the linear equation solver, 
see Listing \ref{lst:bench_2} and results in Figure \ref{fig:bench_2}.

\FloatBarrier
\begin{lstlisting}[caption=Two dimensional diffusion-reaction benchmark with inefficient implementation., label=lst:bench_1, style=Julia]
using OrdinaryDiffEqRosenbrock, Plots, DiffEqDevTools, LinearAlgebra

function f!(du, u, JD, t) #-- RHS: linear diffusion + nonlinear reaction.
   mul!(du, JD, u)  
   du .+= u.^2 .* (1 .- u) .+ exp(t)
   return nothing
end

function J_diff(N, h, D) #-- diffusion matrix
   J = zeros(N^2, N^2)  
   for i = 1:N^2
      J[i,i] = -4
      if ((i-1) >= 1) & (mod(i,N) != 1) J[i,i-1] = 1; end
      if ((i+1) <= N^2) & (mod(i,N) != 0) J[i,i+1] = 1; end
      if (i-N) >= 1 J[i,i-N] = 1; end
      if (i+N) <= N^2 J[i,i+N] = 1; end
   end
   J .*= D / h^2
end

N = 10; final_t = 1.0;  D = 0.1; h = 1 / (N + 1)
JD = J_diff(N, h, D)
u0 = [16*(h * i) * (h * j) * (1 - h * i) * (1 - h * j) for i in 1:N for j in 1:N]

f_ode = ODEFunction(f!)
prob = ODEProblem(f_ode, u0, (0.0, final_t), JD)

ref_sol = solve(prob, Rodas5P(), abstol = 1.0e-14, reltol = 1.0e-14)

numrun = 5
abstols = 1.0 ./ 10.0 .^ (4:11); reltols = 1.0 ./ 10.0 .^ (4:11)
setups = [Dict(:alg=>Rodas4P()),Dict(:alg=>Rodas5P()),Dict(:alg=>Rodas6P()),Dict(:alg=>Rodas4PW(;jac_reuse_gamma_tol=100.0,max_jac_age=100))] 

wp = WorkPrecisionSet(prob,abstols,reltols,setups; save_everystep=true, error_estimate=:l2, appxsol=ref_sol,maxiters=Int(1e6),numruns=numrun)
plot(wp)
\end{lstlisting}

\FloatBarrier
\begin{lstlisting}[caption=Modification of implementation of the two dimensional diffusion-reaction benchmark with sparse matrix operations., label=lst:bench_2, style=Julia]
using OrdinaryDiffEqRosenbrock, Plots, DiffEqDevTools, LinearAlgebra, SparseArrays, LinearSolve

function J_diff(N, h, D) #-- diffusion matrix
   J = spzeros(N^2, N^2)  #-- sparse
   ...
end

N = 60

f_ode = ODEFunction(f!, jac_prototype = (JD .!= 0) )

setups = [Dict(:alg=>Rodas4PW(linsolve = KLUFactorization();jac_reuse_gamma_tol=100.0,max_jac_age=100)), Dict(:alg=>Rodas4P(linsolve = KLUFactorization())),Dict(:alg=>Rodas5P(linsolve = KLUFactorization())),Dict(:alg=>Rodas6P(linsolve = KLUFactorization()))] 

\end{lstlisting}

In both cases, {\tt Rodas4PW} is called with the parameters \verb|max_jac_age=100| and \verb|jac_reuse_gamma_tol=100.0|. This means that the Jacobian matrix should be kept constant for as long as possible. 
The matrix is updated only every 100 time steps or when a step is rejected.
In the first implementation, it is advantageous that {\tt Rodas4PW} requires significantly fewer evaluations of the Jacobian matrix than the ROW methods 
that are not W-methods. Even among the tested W-methods, {\tt Rodas4PW} is the most efficient.
However, when using a sparse Jacobian matrix, the advantage of the W-methods becomes increasingly less significant.

\begin{figure}[t]
 \centering
 \includegraphics[width=0.8\textwidth]{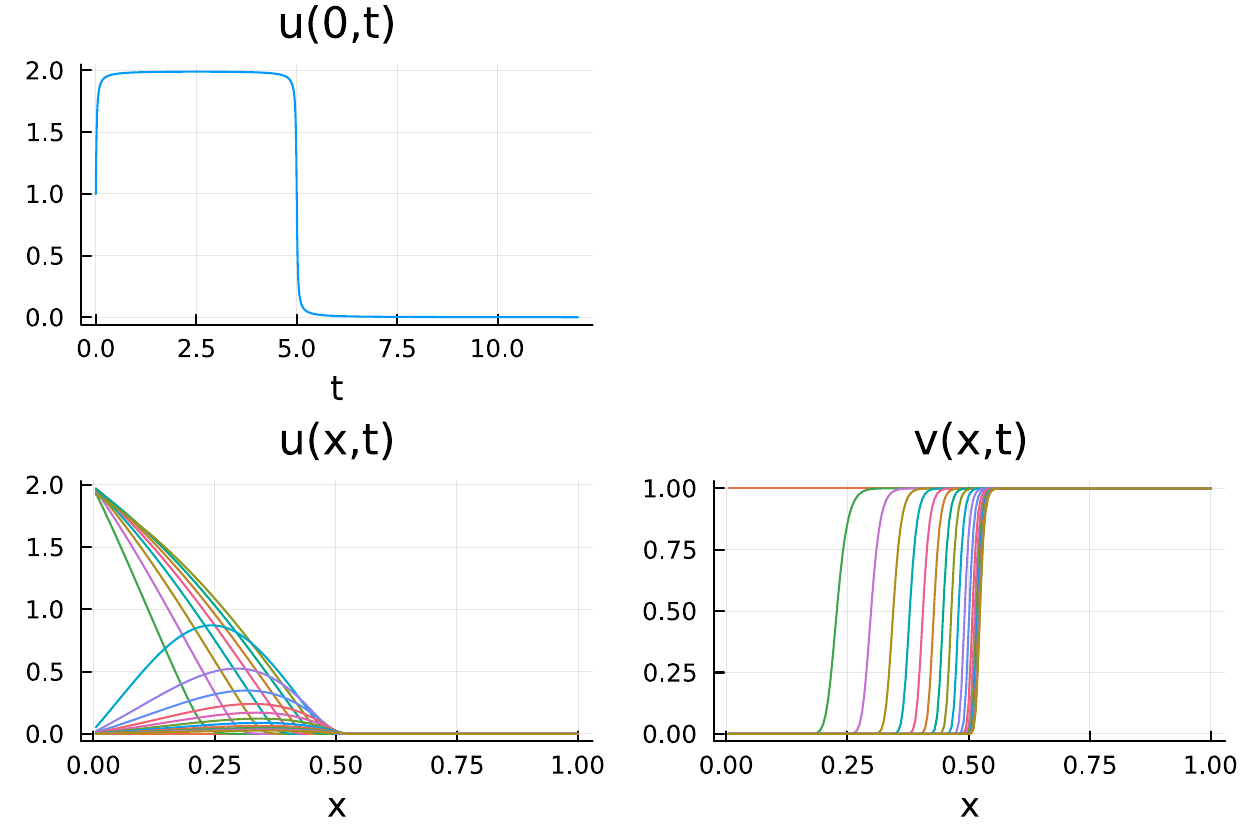} 
 \caption{Boundary condition (\ref{eq:rb}) and solution behavior at different time steps of problem (\ref{eq:prob3a}), (\ref{eq:prob3b}).}\label{fig:hunds}
\end{figure}

The final example is taken from Hundsdorfer \& Verwer \cite{hundsdorfer}. It is the coupled one dimensional advection-diffusion-reaction equation
\begin{eqnarray}
\frac{\partial u}{\partial t} &=& \frac{1}{16} (1-x)^4 \frac{\partial^2 u}{\partial x^2} - \frac{1}{8} (1-x)^3 \frac{\partial u}{\partial x} - 100 \, u \, v \, , \label{eq:prob3a} \\
\frac{\partial v}{\partial t} &=&- 100 \, u \, v \,. \label{eq:prob3b}
\end{eqnarray}
for $x \in (0,1]$ and $t \in (0,12]$ with constant initial values $u(x,0)=0$, $v(x,0)=1$. Since the advection and diffusion terms vanish at $x=1$ the definition 
of a right boundary condition is not necessary.
Instead of the discontinuous left-hand boundary condition for $u$ proposed in \cite{hundsdorfer}, we use the smoothed version:

\begin{equation}
u(0,t) = \frac{2}{\pi} \left( \arctan(50 t) - \arctan(50(t - 5)) \right). \label{eq:rb}
\end{equation}
Beside $u(0,t)$ the solution behavior of $u(x,t)$ and $v(x,t)$ for different times is shown in Figure \ref{fig:hunds}.

\begin{figure}[t]
 \centering
 \includegraphics[width=0.49\textwidth]{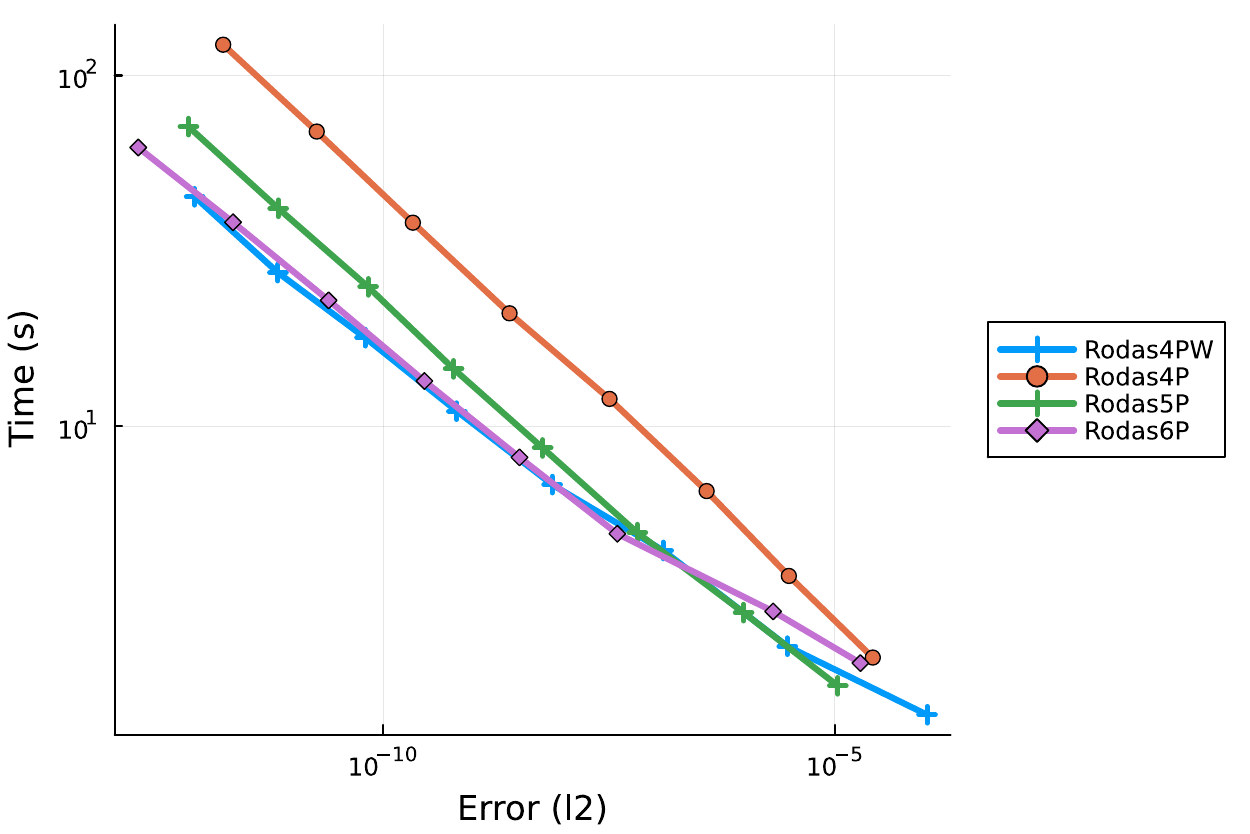} \hfill
 \includegraphics[width=0.49\textwidth]{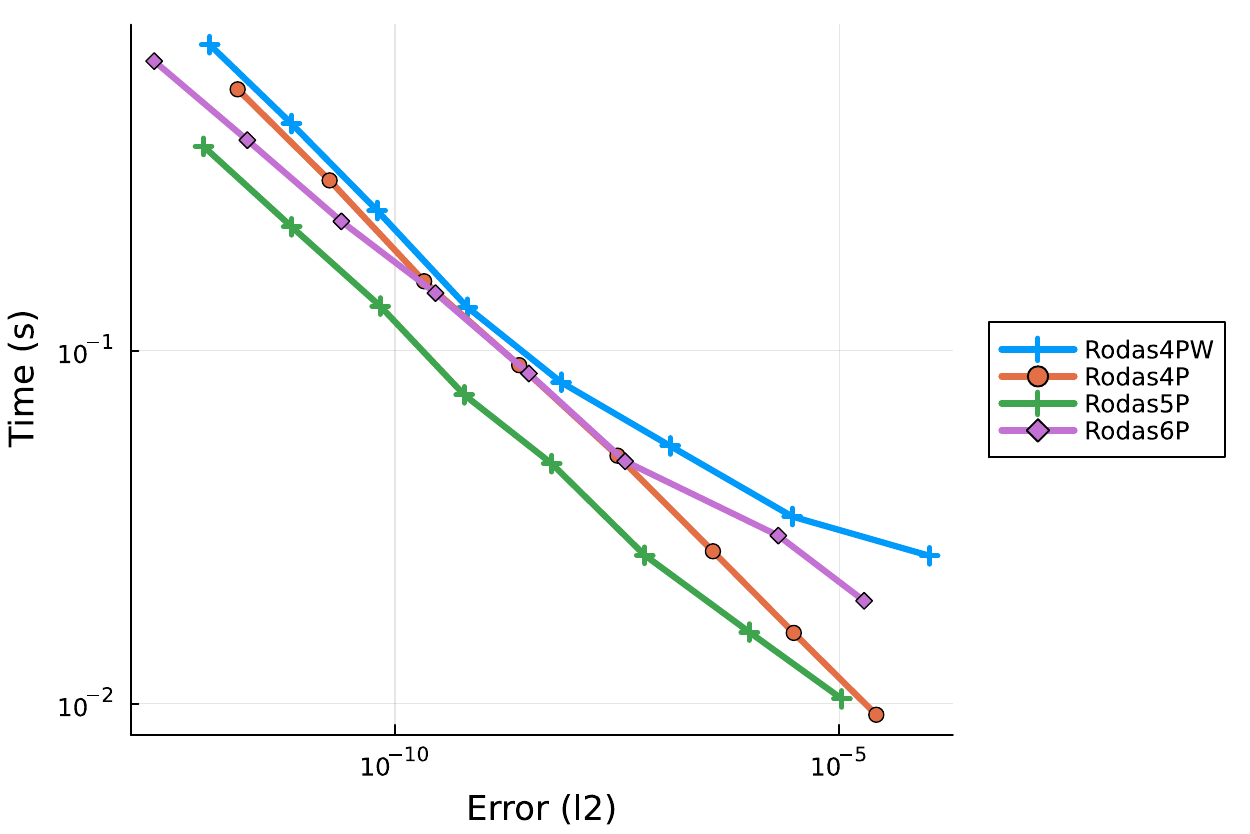}
 \caption{Results of one dimensional advection-diffusion-reaction benchmark with dense matrix (left) and with sparse matrix operations (right).}\label{fig:hundsres}
\end{figure}

As in the previous example, we discretise the spatial derivatives using central differences and create a non-efficient version using dense matrices 
and an efficient version using sparse matrices.
The corresponding Listings and results are shown in Listings \ref{lst:hunds1}, \ref{lst:hunds2} and Figure \ref{fig:hundsres}.

\FloatBarrier
\begin{lstlisting}[caption=One dimensional advection-diffusion-reaction benchmark with inefficient implementation., label=lst:hunds1, style=Julia]
using OrdinaryDiffEqRosenbrock, Plots, DiffEqDevTools, LinearAlgebra

function f!(dy, y, P, t) 
   N, h, k, JD = P
   u = view(y, 1:N); v = view(y, N+1:2*N)
   du = view(dy, 1:N); dv = view(dy, N+1:2*N)
   mul!(du, JD, u)  
   u_0 = 2/pi*(atan(50*t) - atan(50*(t - 5)))
   du[1] += u_0*((1 - h)^4 / (16*h^2) + (1 - h)^3 / (16*h))
   dv .= -k*u.*v  
   du .+= dv  
   return nothing
end

function J_diff(N, h, x) #-- advection, diffusion matrix
   J = zeros(N, N)  
   for i = 1:N
      advec = -(1 - x[i])^3 / (8*2*h) 
      diffus = (1 - x[i])^4 / (16*h^2) 
      J[i,i] = -2*diffus
      if (i-1 >= 1) J[i,i-1] = diffus - advec; end
      if (i+1 <= N) J[i,i+1] = diffus + advec; end
   end
   J
end

N = 200; k = 100; final_t = 12.0; h = 1 / N
x =  LinRange(h, 1.0, N); JD = J_diff(N, h, x)
u0 = [zeros(N); ones(N)]; P = N, h, k, JD

f_ode = ODEFunction(f!)
prob = ODEProblem(f_ode, u0, (0.0, final_t), P)

ref_sol = solve(prob, Rodas5P(), abstol = 1.0e-14, reltol = 1.0e-14)
println(ref_sol.stats)

numrun = 5
abstols = 1.0 ./ 10.0 .^ (4:11); reltols = 1.0 ./ 10.0 .^ (4:11)

setups = [Dict(:alg=>Rodas4PW(;jac_reuse_gamma_tol=100.0,max_jac_age=100)), Dict(:alg=>Rodas4P()), Dict(:alg=>Rodas5P()), Dict(:alg=>Rodas6P())] 
wp = WorkPrecisionSet(prob,abstols,reltols,setups; save_everystep=true, error_estimate=:l2, appxsol=ref_sol,maxiters=Int(1e6),numruns=numrun)
plot(wp)
\end{lstlisting}

\FloatBarrier
\begin{lstlisting}[caption=Modification of implementation of the one dimensional advection-diffusion-reaction benchmark with sparse matrix operations., label=lst:hunds2, style=Julia]
using OrdinaryDiffEqRosenbrock, Plots, DiffEqDevTools, LinearAlgebra, SparseArrays, LinearSolve

function J_diff(N, h, x) #-- advection, diffusion matrix
   J = spzeros(N, N)  #-- sparse
   ...
end
   
Id = sparse(Matrix{Float64}(I, N, N))
J_pattern = [JD Id;Id Id] .!= 0
f_ode = ODEFunction(f!, jac_prototype = J_pattern)

setups = [Dict(:alg=>Rodas4PW(linsolve = KLUFactorization();jac_reuse_gamma_tol=100.0,max_jac_age=100)), Dict(:alg=>Rodas4P(linsolve = KLUFactorization())),Dict(:alg=>Rodas5P(linsolve = KLUFactorization())),Dict(:alg=>Rodas6P(linsolve = KLUFactorization()))] 
\end{lstlisting}

We observe behaviour similar to that in the previous example. The more expensive it is to compute the Jacobian matrix, the more beneficial the W-method becomes. 
However, if the sparse Jacobian matrix is utilised, conventional ROW methods - which require fewer time steps - are preferable.
This is evident from the number of time steps, as well as from the evaluations of the Jacobian matrix and the functional evaluations on the right-hand 
side of (\ref{eq:dae}), which are shown in Table \ref{tab:stat}.
It should be noted that the optimal choice of parameters \verb|jac_reuse_gamma_tol| and \verb|max_jac_age| for this problem has not been examined in detail.
However, no set of parameters could be found that would make {\tt Rodas4PW} operate more efficiently than {\tt Rodas5P} in the case of sparse Jacobian.

\begin{table}
\caption{Solver statistics for problem (\ref{eq:prob3a}), (\ref{eq:prob3b}), calculated with AbsTol=RelTol=$10^{-6.}$. 
NStep = accepted steps, NFail = rejected steps, NJac = Jacobian evaluations,\\ NFcn = function evaluations with dense and sparse Jacobian.}
\label{tab:stat}
\begin{scriptsize}
\begin{center}
\begin{tabular}{|l|r|r|r|r|r|}
  \hline
   & \textbf{NStep} & \textbf{NFail} & \textbf{NJac} & \textbf{NFcn}-dense & \textbf{NFcn}-sparse \\
  \hline
  \textbf{Rodas4P} & 367 &13& 367 &149443 & 4117 \\
  \textbf{Rodas5P} & 279 &14& 279 & 114225 & 3741 \\
  \textbf{Rodas6P} & 259 &19& 259 &109143 & 6579 \\
  \textbf{Rodas4PW}, \verb|max_jac_age=1| & 399 &16& 415 & 170155 & 5812 \\
  \textbf{Rodas4PW}, \verb|max_jac_age=100| & 662 &105& 112 & 51817 & 7465 \\
  \hline
\end{tabular}
\end{center}
\end{scriptsize}
\end{table}


\section{Conclusion}

With {\tt Rodas4PW}, a new W-method of order $p=4$ was derived, which appears to be more efficient than previously known W-methods in many cases. 
However, this paper did not investigate how a general-purpose implementation for exploiting the W-property should be carried out, and this likely 
depends heavily on the specific problem at hand.
If one takes advantage of the fact that, in many cases, the Jacobian matrix is sparse, conventional ROW methods are often the better choice.

%

\newpage
\appendix
\section{Appendix}
\begin{table}[H] 
\caption{Order conditions for W-methods up to order $p=4$ for DAEs assumed that the derivatives of the algebraic equations are exact.}
\label{tab:rodas4pw}
\end{table}
\vspace*{-0.5cm}
\begin{tabular}{lcl}
\hline
1. & $\vcenter{\hbox{\includegraphics[width = 0.075\textwidth]{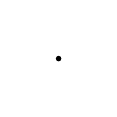}}}$ & $\;\sum b_i= \frac{1}{1}$\\
2. & $\vcenter{\hbox{\includegraphics[width = 0.075\textwidth]{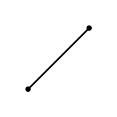}}}$ & $\;\sum b_i\alpha_{ij}= \frac{1}{2}$\\
3. & $\vcenter{\hbox{\includegraphics[width = 0.075\textwidth]{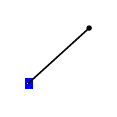}}}$ & $\;\sum b_i\gamma_{ij}= 0$\\
4. & $\vcenter{\hbox{\includegraphics[width = 0.075\textwidth]{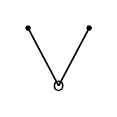}}}$ & $\;\sum b_iw_{ij}\alpha_{jk}\alpha_{jl}= \frac{1}{1}$\\
5. & $\vcenter{\hbox{\includegraphics[width = 0.075\textwidth]{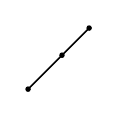}}}$ & $\;\sum b_i\alpha_{ij}\alpha_{jk}= \frac{1}{6}$\\
6. & $\vcenter{\hbox{\includegraphics[width = 0.075\textwidth]{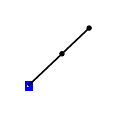}}}$ & $\;\sum b_i\gamma_{ij}\alpha_{jk}= 0$\\
7. & $\vcenter{\hbox{\includegraphics[width = 0.075\textwidth]{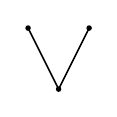}}}$ & $\;\sum b_i\alpha_{ij}\alpha_{ik}= \frac{1}{3}$\\
8. & $\vcenter{\hbox{\includegraphics[width = 0.075\textwidth]{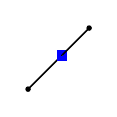}}}$ & $\;\sum b_i\alpha_{ij}\gamma_{jk}= 0$\\
9. & $\vcenter{\hbox{\includegraphics[width = 0.075\textwidth]{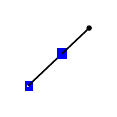}}}$ & $\;\sum b_i\gamma_{ij}\gamma_{jk}= 0$\\
10. & $\vcenter{\hbox{\includegraphics[width = 0.075\textwidth]{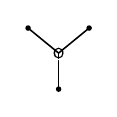}}}$ & $\;\sum b_i\alpha_{ij}w_{jk}\alpha_{kl}\alpha_{km}= \frac{1}{3}$\\
11. & $\vcenter{\hbox{\includegraphics[width = 0.075\textwidth]{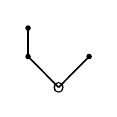}}}$ & $\;\sum b_iw_{ij}\alpha_{jk}\alpha_{jl}\alpha_{lm}= \frac{1}{2}$\\
12. & $\vcenter{\hbox{\includegraphics[width = 0.075\textwidth]{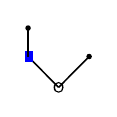}}}$ & $\;\sum b_iw_{ij}\alpha_{jk}\alpha_{jl}\gamma_{lm}= 0$\\
13. & $\vcenter{\hbox{\includegraphics[width = 0.075\textwidth]{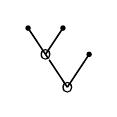}}}$ & $\;\sum b_iw_{ij}\alpha_{jk}\alpha_{jl}w_{lm}\alpha_{mn}\alpha_{mo}= \frac{1}{1}$\\
14. & $\vcenter{\hbox{\includegraphics[width = 0.075\textwidth]{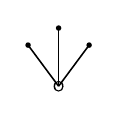}}}$ & $\;\sum b_iw_{ij}\alpha_{jk}\alpha_{jl}\alpha_{jm}= \frac{1}{1}$\\
15. & $\vcenter{\hbox{\includegraphics[width = 0.075\textwidth]{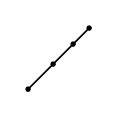}}}$ & $\;\sum b_i\alpha_{ij}\alpha_{jk}\alpha_{kl}= \frac{1}{24}$\\
16. & $\vcenter{\hbox{\includegraphics[width = 0.075\textwidth]{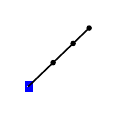}}}$ & $\;\sum b_i\gamma_{ij}\alpha_{jk}\alpha_{kl}= 0$\\
17. & $\vcenter{\hbox{\includegraphics[width = 0.075\textwidth]{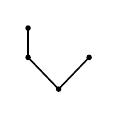}}}$ & $\;\sum b_i\alpha_{ij}\alpha_{ik}\alpha_{kl}= \frac{1}{8}$\\
18. & $\vcenter{\hbox{\includegraphics[width = 0.075\textwidth]{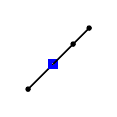}}}$ & $\;\sum b_i\alpha_{ij}\gamma_{jk}\alpha_{kl}= 0$\\
\hline
\end{tabular}
\newpage
\begin{tabular}{lcl}
\hline
19. & $\vcenter{\hbox{\includegraphics[width = 0.075\textwidth]{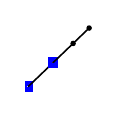}}}$ & $\;\sum b_i\gamma_{ij}\gamma_{jk}\alpha_{kl}= 0$\\
20. & $\vcenter{\hbox{\includegraphics[width = 0.075\textwidth]{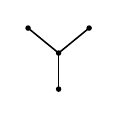}}}$ & $\;\sum b_i\alpha_{ij}\alpha_{jk}\alpha_{jl}= \frac{1}{12}$\\
21. & $\vcenter{\hbox{\includegraphics[width = 0.075\textwidth]{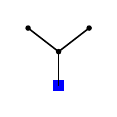}}}$ & $\;\sum b_i\gamma_{ij}\alpha_{jk}\alpha_{jl}= 0$\\
22. & $\vcenter{\hbox{\includegraphics[width = 0.075\textwidth]{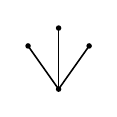}}}$ & $\;\sum b_i\alpha_{ij}\alpha_{ik}\alpha_{il}= \frac{1}{4}$\\
23. & $\vcenter{\hbox{\includegraphics[width = 0.075\textwidth]{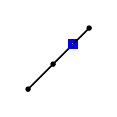}}}$ & $\;\sum b_i\alpha_{ij}\alpha_{jk}\gamma_{kl}= 0$\\
24. & $\vcenter{\hbox{\includegraphics[width = 0.075\textwidth]{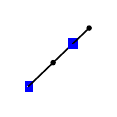}}}$ & $\;\sum b_i\gamma_{ij}\alpha_{jk}\gamma_{kl}= 0$\\
25. & $\vcenter{\hbox{\includegraphics[width = 0.075\textwidth]{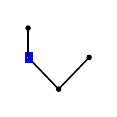}}}$ & $\;\sum b_i\alpha_{ij}\alpha_{ik}\gamma_{kl}= 0$\\
26. & $\vcenter{\hbox{\includegraphics[width = 0.075\textwidth]{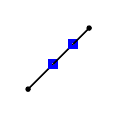}}}$ & $\;\sum b_i\alpha_{ij}\gamma_{jk}\gamma_{kl}= 0$\\
27. & $\vcenter{\hbox{\includegraphics[width = 0.075\textwidth]{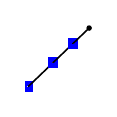}}}$ & $\;\sum b_i\gamma_{ij}\gamma_{jk}\gamma_{kl}= 0$\\
28. & $\vcenter{\hbox{\includegraphics[width = 0.075\textwidth]{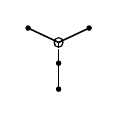}}}$ & $\;\sum b_i\alpha_{ij}\alpha_{jk}w_{kl}\alpha_{lm}\alpha_{ln}= \frac{1}{12}$\\
29. & $\vcenter{\hbox{\includegraphics[width = 0.075\textwidth]{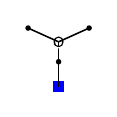}}}$ & $\;\sum b_i\gamma_{ij}\alpha_{jk}w_{kl}\alpha_{lm}\alpha_{ln}= 0$\\
30. & $\vcenter{\hbox{\includegraphics[width = 0.075\textwidth]{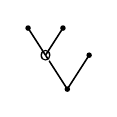}}}$ & $\;\sum b_i\alpha_{ij}\alpha_{ik}w_{kl}\alpha_{lm}\alpha_{ln}= \frac{1}{4}$\\
31. & $\vcenter{\hbox{\includegraphics[width = 0.075\textwidth]{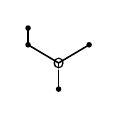}}}$ & $\;\sum b_i\alpha_{ij}w_{jk}\alpha_{kl}\alpha_{km}\alpha_{mn}= \frac{1}{8}$\\
32. & $\vcenter{\hbox{\includegraphics[width = 0.075\textwidth]{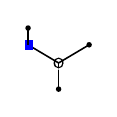}}}$ & $\;\sum b_i\alpha_{ij}w_{jk}\alpha_{kl}\alpha_{km}\gamma_{mn}= 0$\\
33. & $\vcenter{\hbox{\includegraphics[width = 0.075\textwidth]{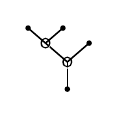}}}$ & $\;\sum b_i\alpha_{ij}w_{jk}\alpha_{kl}\alpha_{km}w_{mn}\alpha_{no}\alpha_{np}= \frac{1}{4}$\\        
34. & $\vcenter{\hbox{\includegraphics[width = 0.075\textwidth]{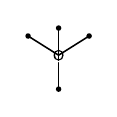}}}$ & $\;\sum b_i\alpha_{ij}w_{jk}\alpha_{kl}\alpha_{km}\alpha_{kn}= \frac{1}{4}$\\
\hline
\end{tabular}

\newpage
\FloatBarrier
\begin{lstlisting}[caption=Coefficients of Rodas4PW., label=lst:coeffs, style=Julia]
alpha = [0.0 0.0 0.0 0.0 0.0 0.0 0.0 0.0 0.0; 
         0.5 0.0 0.0 0.0 0.0 0.0 0.0 0.0 0.0; 
         0.8788271488806836 -0.07252098869793655 0.0 0.0 0.0 0.0 0.0 0.0 0.0; 
         0.1939165447503039 0.21781612476816103 0.13834429354768618 0.0 0.0 0.0 0.0 0.0 0.0; 
         0.07878680633959252 -0.6312867089041129 0.21795734709783576 0.9796662513358836 0.0 0.0 0.0 0.0 0.0; 
        -0.373962028598645 0.1627102517874811 0.33797316465003713 0.15005692742542875 0.46923930452344986 0.0 0.0 0.0 0.0; 
         0.26145229861095143 -0.45561260713981183 -0.4359494907859483 1.06969588654595 -0.4184572253221906 0.3646832990458557 0.0 0.0 0.0; 
         2.036896985284706 -0.09965298531550768 -0.3112689612977009 -1.2549044145263422 0.019552250947243537 0.4625415268297701 0.1468355980778314 0.0 0.0; 
         0.09224774104014404 0.37937607633397097 -0.40989561132505525 -0.4631390612553811 0.15923267359895218 0.44434714536012054 0.54783103624724 0.25 0.0]
gamma = [0.25 0.0 0.0 0.0 0.0 0.0 0.0 0.0 0.0; 
        -0.5 0.25 0.0 0.0 0.0 0.0 0.0 0.0 0.0; 
        -0.6318506050235934 0.319497532555266 0.25 0.0 0.0 0.0 0.0 0.0 0.0; 
        -0.16637199813935943 -0.19027187769208104 -0.13834309258515584 0.25 0.0 0.0 0.0 0.0 0.0; 
        -0.03881495392732915 0.058385689811103814 -0.23411188954598278 -0.35014091687028703 0.25 0.0 0.0 0.0 0.0; 
         0.5527280267671937 -0.6446792849641998 -0.39521691211046756 0.5387809493411774 -0.46924954276022196 0.25 0.0 0.0 0.0; 
        -0.08332750219770624 1.0725563430670497 -0.06638427625294413 -1.5297053034870889 0.18692131378581506 0.27463846464730357 0.25 0.0 0.0; 
        -1.9446492442445618 0.47902906164947867 -0.09862665002735438 0.791765353270961 0.13968042265170866 -0.01819438146964958 0.40099543816940864 0.25 0.0; 
         0.06460887413955387 -0.9420353169841409 0.648812335284155 1.6176460521808094 -0.27791305581935766 -0.7445147175311075 -0.2594646203641316 -0.35713955090575294 0.25]
beta = alpha + gamma
b = beta[9,:]
bd = beta[8,:]
c = [-0.6908483618550232, -6.9322449119244505, 1.512273507237282, 7.634313368493827, 0.2546592489613088, -2.426674669943425, 0.741075199382331, -0.24931715797901102, 0.15676377762716398]
d = [0.49507197684149623, 8.575728497603636, -2.0424811051986396, -10.894508462520267, 0.49657220099816074, 3.2123280759872617, -0.7778179444161315, 0.49863431596037844, 0.43647244474410024]
\end{lstlisting}

\end{document}